\documentclass[a4paper,twoside,
 10pt]{amsart}
\usepackage{amssymb,amsfonts,amsmath,enumerate
}
\newtheorem{theor}{~~~~Theorem}

\newtheorem{prop}{~~~~Proposition}
\newtheorem{cor}{~~~~Corollary}
\newtheorem{lemma}{~~~~Lemma}

\allowdisplaybreaks
\newtheorem{remark}{~~~~Remark}
\newtheorem{defin}{~~~~Definition}
\allowdisplaybreaks
\def\og{\leavevmode\raise.3ex\hbox{$\scriptscriptstyle\langle\!\langle$~}}
\def\fg{\leavevmode\raise.3ex\hbox{~$\!\scriptscriptstyle\,\rangle\!\rangle$}}
\begin{document}
\bibliographystyle{plain}

\title{A note on hyperbolic flows in sub-Riemannian Structures}
\date{\today}
\author{Chengbo Li \address{School of science, Tianjin University, Tianjin, 300072, P.R.China; email: chengboli@gmail.com}}

\begin{abstract}
The \emph{curvature} and the \emph{reduced curvature} are basic
differential invariants of the pair (Hamiltonian system, Lagrange
distribution) on the symplectic manifold. It is shown in
\cite{achamiltonian} that the negativity of the reduced curvature
implies the hyperbolicity of any compact invariant set of the
Hamiltonian flow restricted to a prescribed energy level. We
consider the Hamiltonian flows of the curve of least action of
natural mechanical systems in sub-Riemannian structures with
symmetries.%
%
  We give sufficient conditions for the reduced flows (after reduction of the first
  integrals induced from the symmetries) to be hyperbolic and show
  new examples of Anosov flows.
This result is a generalization of \cite{gmagnetic}  and a partial
 generalization of \cite{mmagnetic} on  magnetic
 flows.


\end{abstract}
\maketitle \markboth {Chengbo Li} {A note on hyperbolic flows  in
sub-Riemannian Structures}
\section{introduction}%

A prime example of Anosov flow is the geodesic flow on a compact
Riemannian manifold with negative sectional curvature
(\cite{ageodesic}). It describes inertial motion of a point particle
confined to the manifold. In this context, magnetic flows, the flows
generated by special forces, were discussed more than 30 years ago
by Anosov and Sinai \cite{ascertain}. They were studied recently by
Gouda \cite{gmagnetic}, Grognet \cite{gflots}, M. and P. Paternain
\cite{panosov} and M. P. Wojtkowski \cite{mmagnetic}. In the last
reference, the potential and the so-called Gaussian thermostats of
external fields were also considered.
.

In the present note, we focus on the hyperbolicity of the flows
associated with a natural mechanical system in a sub-Riemannian
structure with multidimensional symmetries. In this case, the
sub-Riemannian structures are reduced to a Riemannian manifold with
a (vector-valued) magnetic field. We give sufficient conditions for
the reduced Hamiltonian flows (after the reduction of the first
integrals) to be hyperbolic in terms of the Riemannian curvature
tensor and the magnetic field. As a consequence, a class of Anosov
flows are also given.

 In the second section, we formulate the main results of the note. We firstly introduce the notion of a dynamical
Lagrangian distribution and then discuss the reduction after the
first integrals. The key point is that we can construct the
(reduced) curvature maps (forms) for the (reduced) dynamical
Lagrangian distribution based on the work \cite{icdifferential} and
\cite{agfeedback}. The negativity of reduced curvature forms implies
the hyperbolicity of the Hamiltonian flows (\cite{achamiltonian}).
Applying this criteria we give sufficient conditions for the reduced
Hamiltonian flows to be hyperbolic, based on an expression of the
reduced curvature forms via the Riemannian curvature tensor and the
magnetic field.

The last section is devoted to the proofs of the main results. We
apply the similar technique as in \cite{cijacobi} to give the proof
of the expression of the reduced curvature forms and the sufficient
conditions of hyperbolic flows then easily follows.

\section{main results}

\subsection{Dynamical Lagrangian distributions} Let $M$ be
an even dimensional symplectic manifold endowed with a symplectic
form $\sigma$. A \emph{Lagrange distribution}  $\Delta\subset TM$ is
a smooth vector sub-bundle of $TM$ such that each fiber $
\Delta_x=\Delta\cap T_xM,\ x\in M$ is a Lagrangian subspace of the
symplectic space $T_xM$.  Basic examples are cotangent bundles
endowed with the standard symplectic structure and the
\textquotedblleft vertical\textquotedblright distribution:
\begin{equation}\label{exap}
M=T^*N,\ \Pi_x=T_x(T^*_qN),\quad\forall x=(p,q)\in T^*M,p\in T^*_qM,
q\in M.
\end{equation}

Let $h$ be a Hamiltonian function on $M$ and denote by $\vec h$ the
corresponding Hamiltonian vector field: $i_{ \vec h}\sigma=dh$. We
will assume that $ \vec  h$ is a complete vector field without loss
of generality since we will study the dynamics of the Hamiltonian
systems on a compact set. The pair $( \vec  h,\Delta)$ will be said
to be \emph{a dynamical Lagrangian distribution} of the symplectic
manifold $(M,\sigma).$

Dynamical Lagrangian distributions appear naturally in Differential
Geometry, Calculus of Variations and Rational Mechanics. The model
example can be described as follows:

\vskip .1in
 {\bf Example 1}\ On a manifold
$M$ for a given smooth function $L : TM\rightarrow\mathbb R$, which
is convex on each fiber, we consider the following standard problem
of Calculus of Variation with fixed endpoints $q_0$ and $q_1$ and
fixed time $T$:

\begin{eqnarray}\label{action}
&&A(q(\cdot))=\int^T_0 L(q(t),\dot q(t))dt\mapsto {\rm min}\\
&&q(0)=q_0,\quad q(T)=q_1.
\end{eqnarray}
Suppose that the Legendre transform $h : T^*M\rightarrow \mathbb R$
of the function $L$,
\begin{equation}\label{Legendre}
h(p, q) =\max_{X\in T_qM} (p(X)-L(q,X)),\ q\in M, p\in T^*_qM
\end{equation}
 is well defined and smooth on $T^*M$. We will say that the dynamical
Lagrangian distributions $(\vec  h,\Pi)$ is associated with the
problem (2)-(3), where $\Pi$ is as in \eqref{exap}.\quad$\square$

\vskip .1in
 To describe the dynamical property of a dynamical Lagrangian distribution $(\vec  h,\Delta)$, we define the Jacobi
curve (at point $x\in M$) of the pair $( \vec  h,\Delta)$ as
follows:
\begin{equation}\label{Jacur1}
J_x(t):=e^{-t \vec  h}_*\left(\Delta_{e^{t \vec  h}x}\right),
\end{equation}
where $e^{t \vec  h},\ t\in\mathbb R$ denotes the Hamiltonian flow
generated by the vector field $ \vec  h$.

It is clear that the Jacobi curves $J_x(t)$ are curves in the
Lagrange Grassmannian of the symplectic space $T_x^*M$. They are not
arbitrary curves of the Lagrangian Grassmannian but inherit special
features of the pair $( \vec  h,\mathcal D)$. To specify these
features recall that the tangent space $T_\Lambda L(W)$ to the
Lagrangian Grassmannian $L(W)$ of a linear symplectic space $W$
(endowed with a symplectic form $\omega$) at the point $\Lambda$ can
be naturally identified with the space ${\rm Quad}(\Lambda)$ of all
quadratic forms on linear space $\Lambda\subset W$.
Namely, given $\mathfrak V\in T_\Lambda L(W)$ take a curve
$\Lambda(t)\in L(W)$ with $\Lambda(0)=\Lambda$ and
$\dot\Lambda=\mathfrak V$. Given some vector $l\in\Lambda$, take a
curve $\ell(\cdot)$ in $W$ such that $\ell(t)\in \Lambda(t)$ for all
$t$ and $\ell(0)=l$. Define the quadratic form
\begin{equation}
\label{quad} Q_{\mathfrak V}(l)=\omega(l,\frac{d}{dt}\ell(0)).
\end{equation}
Using the fact that the spaces $\Lambda(t)$ are Lagrangian,
it is easy to see that $Q_{\mathfrak V}(l)$ does not depend on the
choice of the curves $\ell(\cdot)$  and $\Lambda(\cdot)$ with the
above properties, but depends only on $\mathfrak V$.
So, we have the linear mapping from $T_\Lambda L(W)$ to the spaces
${\rm Quad}(\Lambda)$, $\mathfrak V\mapsto Q_{\mathfrak V}$.
A simple counting of dimensions shows that this mapping is a
bijection and it defines the required identification. A curve
$\Lambda(\cdot)$ in a Lagrange Grassmannian is called
\emph{regular}, if its velocity  is a nondegenerated quadratic form
at every $\tau$. A curve $\Lambda(\cdot)$ is called \emph{monotone}
(monotonically nondecreasing or monotonically nonincreasing) if the
velocity is sign definite (nonnegative or nonpositive) at any point.
For later convenience, a dynamical Lagrangian distribution is said
to be regular (monotone) if the associated Jacobi curves are regular
(monotone).
%

The group of symplectomorphisms of the ambient space acts naturally
on Lagrangian distribution and Hamiltonian vector fields, therefore
it acts also on dynamical Lagrangian distributions.
 It turns out (\cite{icdifferential}) that one can construct the
canonical bundle of moving frames and the complete system of
symplectic invariants for parametrized curves in Lagrange
Grassmannians satisfying very general assumptions (including
monotone curves as a particular case). The complete system of
symplectic invariants (value at $t=0$) for the Jacobi curve
$J_x(\cdot),x\in M$ is called the curvature maps of $( \vec
h,\Delta)$ and it is the basic differential invariants of the pair
$( \vec h,\Delta)$ w.r.t. the action of symplectic group of $M$. In
this section, we will restrict us to the curvature maps for monotone
regular
 dynamical Lagrangian distribution
  since our goal is to obtain a sufficient
condition for hyperbolicity of the reduced Hamiltonian flows after
the reduction of first integrals, while the reduced dynamical
Lagrangian distributions are monotone regular (see Lemma
\ref{3regular} below). Note also that the curvature maps for regular
curves in Lagrangian Grassmannians are constructed in earlier work
\cite{agfeedback}.

More precisely, let $\mathfrak R_x(t)$ be the curvature map for the
Jacobi curve $J_x(t),x\in M$. Then the linear map $\mathfrak R_x^{(
\vec h,\Delta)}:=\mathfrak R_x(0)=\mathfrak
R_x(t)|_{t=0}:\Delta_x\rightarrow\Delta_x$ is said to be the
\emph{curvature map} (at $x$) of the dynamical Lagrangian
distribution $( \vec h,\Delta)$. It gives a symmetric bilinear forms
(at $x$)
$$r^{( \vec  h,\Delta)}_x(v,w):=\sigma(R_x^{( \vec
h,\Delta)}w,[ \vec  h,V]),\quad v,w\in\Delta_x$$
 where $V$ is a smooth section of the
sub-bundle $\Delta$ with $V(x)=v$. The corresponding quadratic form
will be called the \emph{curvature form} of the dynamical Lagrangian
distribution $( \vec  h,\Delta)$.

\vskip .1in
 {\bf Example 2}
(Natural mechanical system) In Example 1, let
$$M =
R^n ,\Pi_{(p,q)} = (R^n, 0), L(q,X) = \frac{1}{2}|X|^2 - W(q)$$ (in
this case the function $A(q(t),\dot q(t))$  is the Action functional
of the natural mechanical system with potential energy $W(q)$). Then
the curvature forms can be written as follows:
\begin{equation}
r^{( \vec
h,\Pi)}_{(p,q)}(\partial_{p_i},\partial_{p_j})=\frac{\partial^2W}{\partial
q_i\partial q_j}(q),\quad\forall 1\leq i,j\leq n.
\end{equation}
 In other words, in
this case the curvature forms are naturally identified with the
Hessian of the potential $W$.\quad$\square$

\vskip .1in
 {\bf  Example 3} (Riemannian manifold) Let $(M,g)$ be a Riemannian
manifold. Let $L(q,X) = \frac{1}{2}g(X,X)$. The inner product
$g(\cdot, \cdot)$ defines the canonical isomorphism between $T_qM$
and $T^*_qM$. For any $q\in M$ and $p\in T^*_qM$ we will denote by
$p^h$ the image of $p$ under this isomorphism, namely, the vector
$p^h\in T_qM$, satisfying
\begin{equation}\label{iso}
p(\cdot) = g(p^h, \cdot)
\end{equation}
Since the fibers $T_q^*M$ are linear spaces, one can identify
$\Pi_\lambda (= T_\lambda T^*_qM)$ with $T_{\pi(\lambda)}^*M$, i.e.
the operation $p^h$ is defined also on each $p\in\Pi_\lambda$ with
values in $T_{\pi(\lambda)}M$. For any given $\lambda=(p,q)\in
T^*M,p\in M,p\in T^*_qM$, it turns out (\cite{agfeedback}) that
\begin{equation}\label{Riemannian} (\mathfrak
R^{(\vec  h,\Pi)}_\lambda v)^h=R^\nabla(p^h,v^h)p^h,\quad
v\in\Pi_\lambda,
\end{equation}
 where $R^\nabla$ is the
Riemannian curvature tensor of the metric $g$.\quad$\square$ 
\medskip

{\bf Example 4} (Natural mechanical system on a Riemannian manifold)
We add the potential in the action functional in the previous
example, i.e. $L(q,X) = \frac{1}{2}g(X,X)-W(q)$ . 
Then the curvature maps satisfies
\begin{equation}\label{RiePo}
(\mathfrak R^{(\vec  h,\Pi)}_\lambda
v)^h=R^\nabla(p^h,v^h)p^h+\nabla_{v^h}(\nabla W)(q),\quad
v\in\Pi_\lambda,
\end{equation}
where $\nabla W$ is the gradient of $W$ w.r.t. the Riemannian metric
$g$.\quad$\square$

\subsection{Reduced curvature forms and hyperbolicity} Assume that the
dynamical Lagrange distribution $( \vec  h,\Delta)$ have arbitrary
$s$ first integrals $g_1,...,g_s$ in involution with the Hamiltonian
$h$, i.e. $s$ functions on $M$ such that
$$\{h,g_i\}=0,\quad \{g_i,g_j\}=0,\quad \forall 1\leq i,j\leq s,$$
where $\{,\}$ is the Poisson bracket. This problem appear naturally
in the framework of mechanical systems and variational problems with
symmetries. Let $\mathcal G=(g_1,...,g_s)$ and let
\begin{equation}\label{redu1}
\Delta^{\mathcal G}_x=(\cap_{i=1}^s {\rm ker}\
d_xg_i)\cap\Delta_x+{\rm span}\{\vec g_1(x),...,\vec g_s(x)\}
\end{equation}

Clearly, the distribution $\Delta^{\mathcal G}=\{\Delta^{\mathcal
G}_x,x\in M\}$ is a Lagrangian distribution. Hence, we get a reduced
dynamical Lagrangian distribution $( \vec h,\Delta^{\mathcal G})$
after the reduction by first integrals $\mathcal G$. Its curvature
maps (forms) will be called \emph{the reduced curvature maps
(forms)} after the reduction by first integrals $\mathcal G$.
%
%

\vskip .1in
 {\bf Example 5} Assume that we have one first integral $g$ of $h$
such that the Hamiltonian vector field $\vec  g$ preserves the
distribution $\Delta$, i.e. $(e^{t\vec g})_*\Delta=\Delta$. Fixing
some value $c$ of $g$, one can define (at least locally) the
following quotient manifold: $M_{g,c} = g^{-1}(c)/{\mathcal C}$ ,
where $\mathcal C$ is the line foliation of the integral curves of
the vector field $\vec  g$. The manifold $M_{g,c}$ naturally
inherits a symplectic form from the original symplectic structure
$(M,\sigma)$. Furthermore, if we denote by $\Phi :
g^{-1}(c)\rightarrow M_{g,c}$ the canonical projection on the
quotient set, the vector field $\Phi_*(\vec  h)$ is well defined
Hamiltonian vector field on $M_{g,c}$ due to the fact that the
vector fields $\vec  h$ and $\vec g$ commute. For simplicity, we
still denote $\Phi_*(\vec  h)$ by $\vec h$. Actually we have simply
described the standard reduction of the Hamiltonian systems on the
level set of the first integrals in Mechanics (see e.g.
\cite{amfoundations}). In this way, to any dynamical Lagrangian
distribution $(\vec h,\Delta)$ on $M$ one can associate the
dynamical Lagrangian distribution $(\vec h,\Phi_*\Delta)$ on the
symplectic manifold $M_{g,c}$ of smaller dimension.\quad$\square$
%
%

It is well known that the geodesic flows on a compact Riemannian
manifold with negative sectional curvature is Anosov
(\cite{ageodesic}). On the other hand, the reduced curvature maps
(forms) of the dynamical Lagrangian distributions associated with
the geodesic problem on a Riemannian manifold are naturally
identified with the sectional curvature tensor (See Example 6
below). Hence, we could roughly formulate the result of Anosov as
follows: negativity of the reduced curvature forms implies the
hyperbolicity of the geodesic flows. To go further from this
viewpoint, one can obtain a natural generalization
(\cite{achamiltonian})  in the framework of dynamical Lagrangian
distributions. It will serve as a criteria in the study of the
hyperbolic flows in sub-Riemannian structures.
\begin{defin}
Let $e^{tX},\ t\in\mathbb R$ be the flow generated by the vector
field $X$ on a manifold $P$. A compact invariant set $A\subset P$ of
the flow $e^{tX}$ is called a hyperbolic set if there exists a
Riemannian structure in a neighborhood of $A$, a positive constant
$\delta$, and a splitting: $T_zP=E_z^+\oplus E_z^-\oplus\mathbb
RX(z),\ z\in A$ such that $X(z)\neq 0$ and
\begin{enumerate}
\item $e^{tX}_*E^+_z=E^+_{e^{tX}z},\ e^{tX}_*E^-_z=E^-_{e^{tX}z},$
\item $\|e^{tX}_*\zeta^+\|\geq e^{\delta t}\|\zeta^+\|,\
\forall t>0, \forall \zeta^+\in E^+_z,$
\item $\|e^{tX}_*\zeta^-\|\leq e^{-\delta t}\|\zeta^-\|,\
\forall t>0, \forall \zeta^-\in E^-_z.$
\end{enumerate}
If the entire manifold $P$ is a hyperbolic set, then the flow
$e^{tX}$ is called a flow of Anosov type.
\end{defin}

\begin{theor}\label{main1}
 Let $c=(c_0,c_1,...,c_s)$ be constants. Let
$S$ be a compact invariant set of the flow $e^{t\vec  h}$ contained
in a fixed level of $h^{-1}(c_0)\cap_{i=1}^s g_i^{-1}(c_i)$ and $
\vec  h(x), \vec  g_i(x)\notin \mathcal D_x,\forall x\in
S,i=1,...,s$. If the reduced curvature form $r_x^{(\mathcal G,h)}$
of the dynamical Lagrangian distribution $( \vec  h,\mathcal D)$
(after the reduction by first integrals $(\mathcal G,h$)) is
negative at every point $x$ of $S$, then $S$ is a hyperbolic set of
the flow $e^{t\vec h}|_{h^{-1}(c_0)\cap_{i=1}^s g_i^{-1}(c_i)}$.

\end{theor}%

\subsection{Descriptions of main results}
We now specialize to the study of a natural mechanical system on a
sub-Riemannian manifold with symmetries.

Let $M$ be a connected smooth manifold. A distribution $\mathcal D$
on $M$ is a sub-bundle of the tangent bundle $TM$. It is said to be
\emph{completely nonholonomic} if any local frame $\{X_i: 1\leq
i\leq n\}$ for $\mathcal D$, together with all its iterated Lie
brackets $[X_i, X_j], [X_i, [X_j, X_k]]$, ... , spans the tangent
bundle $TM$.
A Lipschitzian curve $\gamma:[0,T]\longrightarrow M$ is said to be
admissible if $\dot\gamma(t)\in\mathcal D_{\gamma(t)}$ for a.e.
$t\in [0,T]$. From the Rashevskii-Chow theorem (\cite{ascontrol}) it
follows that there is an admissible curve joining any two points of
$M$. \emph{A sub-Riemannian metric} is a smoothly varying positive
definite inner product $\left\langle\cdot,\cdot\right\rangle$ on
$\mathcal D$. In particular, when $\mathcal D$ is equal to the
tangent bundle, $\left\langle\cdot,\cdot\right\rangle$ gives a
Riemannian metric.

\emph{A sub-Riemannian structure}, denoted by the triple $(M,
\mathcal D, \left\langle\cdot,\cdot\right\rangle)$, is a smooth
$n$-dimensional connected manifold $M$ equipped with a
sub-Riemannian metric $\left\langle\cdot,\cdot\right\rangle$ on a
completely nonholonomic distribution $\mathcal D$. In this case, we
call the manifold $M$ \emph{a sub-Riemannian manifold}. In the
present note we consider sub-Riemannian metrics
$\left\langle\cdot,\cdot\right\rangle$ on  distribution $\mathcal D$
of corank $s$, having $s$ transversal infinitesimal symmetries, i.e.
$s$ vector fields  $X_1,\ldots, X_s$ on $M$ such that
\begin{equation}\label{sym}
e^{tX_i}_*\mathcal{D}=\mathcal{D}\ ,\ (e^{tX_i})^*\left\langle\cdot,
\cdot\right\rangle=\left\langle\cdot, \cdot\right\rangle,\quad 1\leq
i\leq s,
\end{equation}
 and $TM=\mathcal D\oplus \text{span}\{X_1,\ldots,X_s\}$.
Suppose further that the symmetries $\{X_i:\ 1\leq i\leq s\}$
 are commutative (see Remark 1 for
noncommutative case), i.e.
\begin{equation}
\label{Xcom} [X_i,X_j]=0,\quad \forall 1\leq i,j\leq s.
\end{equation}

\vskip .1in
 We consider the natural mechanical system on a
sub-Riemannian manifold (ASR):
\begin{eqnarray}\label{actSub}
&&A(\gamma(\cdot))=\int^T_0(\frac{1}{2}\|\dot \gamma\|^2-W(\gamma))dt\mapsto {\rm min}\\
&&\gamma(\cdot)\ {\rm is\ admissible},\quad \gamma(0)=q_0,\quad
\gamma(T)=q_1.
\end{eqnarray}%
where $\|\cdot\|$ is the norm w.r.t. the metric
$\left\langle\cdot,\cdot\right\rangle$.  We assume further that the
potential $W$ in \eqref{actSub} is constant along the integral
curves of any $X_i$, or, equivalently,
\begin{equation}\label{symP}
X_i(W)=0,\ i=1,...,s.
\end{equation}

It is more convenient to regard it as an optimal control problem and
its extremals can be described by the Pontryagin Maximum Principle
of Optimal Control Theory (\cite{pbgmthe}). There are two different
types of extremals: abnormal and normal, according to vanishing or
nonvanishing of Lagrange multiplier near the functional,
respectively. The minimizers of the problem are the projections of
either normal extremals or abnormal extremals.

In the present note we will focus on normal extremals only. To
describe them let us introduce some notations.
Let $T^*M$ be the cotangent bundle of $M$ and $\sigma$ be the
canonical symplectic form on $T^*M$, i.e. $\sigma=-d\varsigma$,
where $\varsigma$ is the tautological (Liouville) 1-form on $T^*M$.
Let
\begin{equation}\label{h}
h(p,q)=\max_{u\in\mathcal{D}}(p\cdot
u-\frac{1}{2}\|u\|^2+W(q))=\frac{1}{2}\|p|_{\mathcal{D}_q}\|^2+W(q),\
q\in M,\ p\in T^*_qM,
\end{equation}
where $p|_{\mathcal{D}_q}$ is the restriction of the linear
functional $p$ to $\mathcal{D}_q$ and the norm
$\|p|_{\mathcal{D}_q}\|$ is defined w.r.t. the Euclidean structure
on $\mathcal D_q.$  It is well defined and smooth in the open set
$O=T^*M\backslash \mathcal D^\perp$, where $\mathcal D^\perp$ is the
annihilator of $\mathcal D$, that is,
\begin{equation}\label{3Dor}
\mathcal D^\perp=\{(p,q)\in T^*M: p(v)=0\,\,\forall v\in \mathcal
D_q\}.
\end{equation}
\medskip
 For any vector field $X_i$ define the \textquotedblleft
quasiimpluses\textquotedblright \ $u_i: T^*M\rightarrow \mathbb R$
by
$$\quad u_i(p,q)=p(X_i(q)),\
q\in T^*_qM, q\in M,\ \forall 1\leq i\leq s.$$ Let $h$ be the
sub-Riemannian Hamiltonian as in \eqref{h}. Then it follows from
\eqref{sym} and \eqref{symP} that
\begin{equation}\label{3symm3}
\{h, u_i\}=0, \forall 1\leq i\leq s.
\end{equation}
and from \eqref{Xcom} it follows that
\begin{equation}\label{Poisson}
\{u_i, u_j\}=0,\ \forall 1\leq i,j\leq s,
\end{equation} where $\{\
,\ \}$ is the Poisson bracket. In other words, $u_i (1\leq i\leq s)$
are first integrals in involution of the Hamiltonian system
$e^{t\vec h}$.

As before, let $\Pi$ be the \textquotedblleft
vertical\textquotedblright distribution, i.e. $\Pi_\lambda=T_\lambda
T_{\pi(\lambda)}^*M$, where $\pi:T^*M\to M$ is the canonical
projection.
 Now we can apply the reduction  to the dynamical
Lagrangian distributions $( \vec  h,\Pi)$
 after the first integrals
$u_i (1\leq i\leq s)$ (c.f. Example 5). For this, fix constants
$c_0,c_1,...,c_s$, where $c_0>0$ is sufficient large. Take a common
level set
$$\mathcal H_c:=\{h=c_0\}\cap\{u_i=c_i, 1\leq i\leq s\}.$$
 Then
 $$W^{c}_{\lambda}=T_{\lambda}\mathcal
H_{c}/{\rm span}\{\vec  h(\lambda), \vec u_i(\lambda),\ 1\leq i\leq
s\}$$
 is a linear symplectic space with the symplectic
form $\sigma^{c}$ naturally inherited from the symplectic form
$\sigma$. Moreover,
$$\Pi_{\lambda}^{c}=(T_{\lambda}\mathcal H_{c}\cap\Pi_{\lambda})/{\rm span}\{\vec  h(\lambda),
\vec  u_i(\lambda),\ 1\leq i\leq s\}$$ is a Lagrangian subspace in
$W^{c}_{\lambda}$. Hence, we get the reduced dynamical Lagrangian
distribution $( \vec  h,\Pi^{c})$ in the linear symplectic space
$W^{c}$.

\begin{remark}
The reduction procedure above also applies for the case that the
symmetries $\{X_i:\ 1\leq i\leq s\}$ are not commutative but still
satisfy that $g=\rm{span}_{\mathbb R}\{X_i,\ 1\leq i\leq s\}$ is a
Lie algebra and the derived Lie algebra $g^2=[g, g]$ is a proper Lie
subalgebra. 
Indeed, take a basis of $g^2:(g_1,..., g_k)$ and complete it to a
basis of $g: (g_1,...,g_k, g_{k+1},...,g_s)$. Then if we select a
level set $c=(c_1,..., c_s)$ such that
\begin{equation}\label{noncom}
c_i=0,\ 1\leq i\leq k,
\end{equation}
then one can see that $g$ is commutative on this level set.
Therefore, it reduces to commutative case and then the Poission
reduction can be applied. Actually, it is equivalent to  considering
a sub-Riemannian structure on the manifold obtained by reduction of
the original one by $g_1,..., g_k$ on which the symmetries consist
of a commutative Lie algebra $g/g^2$.
\end{remark}

%
%

The (reduced) curvature maps (forms) of $(\vec h,\Pi^c)$ is
naturally related to the ambient sub-Riemannian structures (with
symmetries), while the later can be reduced to a Riemannian manifold
equipped with a $\mathbb R^s$-valued magnetic field.
 Denote by $\widetilde M$ the quotient of $M$ by the
leaves of the integral manifold of the involutive distribution
spanned by $X_1,\ldots X_s$ and denote the factorization map by
${\rm pr}: M\rightarrow\widetilde M$. Then $\widetilde M$ is (at
least locally) a Riemannian manifold equipped with the Riemannian
metric $g$ induced from the sub-Riemannian metric. Furthermore, let
$\omega=(\omega_i)_{1\leq i\leq s}$ be the $\mathbb R^s$-valued
1-form defined by $\omega_i|_{\mathcal D}=0$ and
$\omega_i(X_j)=\delta_{ij},\ \forall 1\leq i,j\leq s.$ Then
$d\omega=(d\omega_i)_{1\leq i\leq s}$ induces a $\mathbb R^s$-valued
2-form on $\widetilde M$ (still denoted by $d\omega=(d\omega_i)$)
and one can define a $\mathbb R^s$-valued tensor $J=(J_i(\tilde
q),\tilde q\in\widetilde M)$ of type $(1, 1)$ on $\widetilde M$
satisfying
$$g_{\tilde q}(J_i(\tilde q)v, w)=d\omega_i(\tilde q)(v, w),\ v,w\in T_{\tilde q}\widetilde M, \tilde q\in\widetilde M,\ \forall 1\leq i\leq s.$$

\medskip

Let $\Xi^{c}$ be the $s$-foliation such that its leaves are integral
curves of $\{\vec  u_i, 1\leq i\leq s\}$. Let ${\rm PR}^{c}:T^*M\to
T^*M/\Xi^{c}$ be the canonical projection to the quotient manifold.
Now we show that the quotient manifold $N^{c}=\{u_i=c_i, 1\leq i\leq
s\}/\Xi^{c}$ can be naturally
identified with $T^*\widetilde M$. 
Indeed, a point $\tilde \lambda$ in $\{u_i=c_i,1\leq i\leq
s\}/\Xi^{c}$ can be identified with a leaf
$(\rm{PR}^{c})^{-1}(\tilde\lambda)$ of $\Xi^{c}$ which has a form
$$((e^{-\sum_{i=1}^s t_iX_i})^*p,e^{\sum_{i=1}^s t_iX_i}q),$$ where
$\lambda=(p,q)\in ({\rm PR}^{c})^{-1}(\tilde\lambda)$, $q\in M$ and
$p\in T_q^* M$. On the other hand, any element in $T^* \widetilde M$
can be identified with a one-parametric family of pairs
$((e^{-\sum_{i=1}^s t_iX_i})^*(p|_\mathcal D),e^{\sum_{i=1}^s
t_iX_i}q)$.
The mapping $I^{c}:\{u_i=c_i,1\leq i\leq s\}/\Xi^{c}\to
T^*\widetilde M $ defined by
$$I^{c}:(e^{-\sum_{i=1}^s t_iX_i})^*p,e^{\sum_{i=1}^s
t_iX_i}q)\mapsto (e^{-\sum_{i=1}^s t_iX_i})^*(p|_{\mathcal
D}),e^{\sum_{i=1}^s t_iX_i}q)$$ is one-to-one
($p(X_i)=\hbox{const.}$ is already prescribed) and it defines the
required identification.

Before the statement of the main result of the note, let us
introduce some notations. Let $D^\bot$ be as in \eqref{3Dor}. Denote
$\mathcal D_q^\bot=\mathcal D^\bot\cap T^*_qM.$ Then one has the
following series of natural identifications:
\begin{equation}\label{3ident1}
\Pi^{c}_\lambda\sim T^*_qM/\mathcal{D}_q^\bot\sim
\mathcal{D}^*_q\stackrel{\left\langle\cdot,
\cdot\right\rangle}{\sim}\mathcal{D}_q\sim T_{{\rm pr}(q)}\widetilde
M£¬
\end{equation}
where $\mathcal D_q^*\subseteq T^*_qM$ is the dual space of
$\mathcal D_q$. Given $ v\in T_\lambda T^*_{q}M$ ($\sim T^*_q M$),
where $q=\pi(\lambda)$, we can assign a unique vector $v^h\in
T_{{\rm pr}(q)}\widetilde M$ to its equivalence class in
$T_q^*M/\mathcal D_q^\bot$ by using the identifications
\eqref{3ident1}.
  Conversely, to any $X\in T_{\hbox{pr}(q)}\widetilde M$ one can assign an equivalence class of $T_\lambda(T^*
_qM)/\mathcal D_q^\bot$. Denote by $X^v\in T_\lambda T^*_{q}M$ the
unique representative of this equivalence class such that
$du_i(X^v)=0,\,\forall 1\leq i\leq s$.

For simplicity, we henceforth denote: $J^c=\sum_{i=1}^sc_iJ_i$.

%
\begin{theor}\label{3main1}
The curvature forms $r_\lambda^{c}$ of the dynamical Lagrangian
distribution $(\vec h,\Pi^c)$ is expressed as follows. For any
$v\in\Pi^{c}_\lambda$,
\begin{eqnarray*}
r^{c}_\lambda(v)&=&g(R^{\nabla}(p^h,v^h)p^h, v^h)+ g(\nabla J^c(p^h,
v^h),v^h)+\frac{1}{4}g(J^cv^h,J^cv^h)\\
&+&\frac{3}{8(c_0+W)}\left(g(J^cp^h,v^h)\right)^2+\frac{3}{2(c_0+W)}g(v^h,\nabla
W)g(J^cp^h,v^h)\\
&+&\frac{3}{2(c_0+W)}\left(g(v^h,\nabla W)\right)^2+{\rm Hess}\
W(v^h,v^h).
\end{eqnarray*}
\end{theor}

\medskip
It follows from relations \eqref{symP},\eqref{3symm3} and
\eqref{Poisson} that $e^{t\vec h}$ induces a (reduced) Hamiltonian
flow $\Phi_t$ on $N^c$, where $N^{c}=\{u_i=c_i, 1\leq i\leq
s\}/\Xi^{c}$, as before. 
Then following theorem is a direct consequence Theorem \ref{main1}.
\begin{theor}\label{3main2}%
 Assume that $K^c\subset
N^c$ is a compact invariant set of
the flow $\Phi_t$ on $N^c$. 
 If
the curvature form $r_\lambda^{c}$ is negative at every point of
$K^c$, then $K^c$ is a hyperbolic set of the flow $\Phi_t$ on $N^c$.
\end{theor}%
As mentioned, the manifold $N^c$ is naturally identified $T^*M$. Now
denote by $S_1\widetilde M$ the unit tangent bundle. Combining
the previous theorem with Theorem \ref{3main1},  
we get the
following

\begin{theor}\label{3main3}
Assume that the reduced Riemannian manifold $(\widetilde M, g)$ is
compact and has sectional curvature bounded  from above by $k_{\rm
max}$. If the constants $c_0,c_1,...,c_s$ satisfy
\begin{eqnarray*}
&&\max_{v, w\in S_1\widetilde M, v\perp w}g(v,\nabla
J^c(w;v))+\frac{1}{4}g(J^cv,
J^cv)+\frac{3}{8(c_0+W)}g(w,J^cv)g(w,J^cv)\\
&&+\frac{3}{2(c_0+W)}g(v,\nabla W)g(J^cw,v)+3\left(\frac{\|\nabla
W\|}{2(c_0+W)}\right)^2 +\frac{\|{\rm Hess}\ W\|}{2(c_0+W)}<-k_{\rm
max},
\end{eqnarray*}
then the flow $\Phi_t$
is an Anosov flow.
\end{theor}

\begin{cor}
Pure potential flows, i.e. $J_i=0 (1\leq i\leq s)$,
\begin{equation}\label{PurePo}
\max_{\tilde q\in\widetilde M}\left(3\left(\frac{\|\nabla
W\|}{2(c_0+W)}\right)^2 +\frac{\|{\rm Hess}\
W\|}{2(c_0+W)}\right)<-k_{\rm max}.
\end{equation}
\end{cor}
\begin{cor}
Pure magnetic flows, i.e. $W=0$,
\begin{equation}\label{PureMa}
\max_{v, w\in S_1\widetilde M, v\perp w}{cg(v,\nabla
J(w;v))+c^2g(Jv, Jv))}<-k_{\rm{max}}.
\end{equation}
\end{cor}
\begin{remark}
The left-hand side of the inequality in Theorem \ref{3main3} is
always positive, because the second term inside the $\text{max}$ is
positive and the first term can be made nonnegative, if necessary,
by changing the sign of $w$. Hence, Theorem \ref{3main3} makes sense
only if $k_{\rm {max}}<0$.
\end{remark}
The flow $\Phi_t$ on $N^{c}$ can be considered as a perturbation of
the Riemannian geodesic flow: the flow $\Phi_t$ on $N^c$ remains to
be an Anosov flow for sufficient small constants $c_i (1\leq i\leq
s)$ and for proper potential function $W$ (sufficient small norm of
$W$ and its derivatives). When $s=1$, it coincides with Theorem 4.1
(the case of Gaussian thermostats of external fields $E=0$ there) in
\cite{mmagnetic}; when $s=1$ and $W=0$, the condition of Anosov
magnetic flows \eqref{PureMa} coincides with the main results in
\cite{gmagnetic}.

\section{proof of the main results}
The rest of the note is devoted to the proof of Theorem
\ref{3main1}.

\subsection{Reduced curvature maps}
As before, fix constants $c_0,c_1,...,c_s$, where $c_0>0$ is
sufficient large. Let $J^{c}_{\lambda}(t)$ be the Jacobi curves
associated with the reduced dynamical Lagrangian distribution $(\vec
h,\Pi^c)$, namely,
\begin{equation}\label{redjac}
J^{c}_{\lambda}(t):=e^{-t\vec h}_*\Pi^{c}_{e^{t\vec h}\lambda}.
\end{equation}

\begin{lemma}\label{3regular}
The reduced Jacobi curve $J^{c}_{\lambda}(\cdot)$ is a regular
monotone nondecreasing curve in Lagrange Grassmannian
$L(W^{c}_{\lambda})$.
\end{lemma}
\begin{proof}

First note that if $\bar\lambda=e^{\bar t\vec  h}\lambda$ and $\phi:
W^{c}_\lambda\rightarrow W^{c}_{\bar\lambda}$ is a symplectic
transformation induced in the natural way by a linear mapping
$e^{t\vec  h}_*:T_\lambda \mathcal H_{c}\rightarrow
T_{\bar\lambda}\mathcal H_{c}$, where, as before,
$$\mathcal H_c:=\{h=c_0\}\cap\{u_i=c_i, 1\leq i\leq s\}.$$
Then by (\ref{redjac}) we have
\begin{equation}
\label{3diffpoint}
J^{c}_{\bar\lambda}(t)=\phi\bigl(J^{c}_
\lambda(t-\bar t)\bigr).
\end{equation}
 Further, it turns out (see, for example, \cite[Proposition 1]{azgeometry1}) that the velocity of the Jacobi curve $J_\lambda^{c}(\cdot)$ at $t = 0$ is equal to the restriction of the Hessian of $h$ to the tangent space to $\Pi_{\lambda}^{c}$ at the point $\lambda$.
  This together with the relation \eqref{3diffpoint} and the construction of $\Pi_{\lambda}^{c}$ implies easily that $J^{c}_{\lambda}(t)$ is a regular
  monotone nondecreasing curve.
\end{proof}

\begin{theor}
Let $\Lambda(\cdot)$ be a regular curve in the Lagrange Grassmannian
$L(G)$ of a $2n$-dimensional linear symplectic space $G$. Then there
exists a moving Darboux frame $(E(t), F(t))$ of $G$:
$$E(t)=(e_1(t),...,e_{n}(t)),\ F(t)=(f_1(t),...,f_{n}(t))$$
such that $\Lambda(t)=\rm{span}\{E(t)\}$ and there exists a
one-parametric family of linear self-adjoint operators $\mathfrak
R(t): \Lambda(t)\rightarrow \Lambda(t)$ satisfying
\begin{equation}
\label{3structeq}
\begin{cases}
 E^{\prime}(t)=F(t),\\
 F^{\prime}(t)=-R(t)E(t).
\end{cases}
\end{equation}
The moving frame $(E(t), F(t))$ is a called \emph{a normal moving
frame} of $\Lambda(t)$ and the linear operator $\mathfrak R(t)$ is
called \emph{the curvature map} of $\Lambda(t)$. A moving frame
$(\widetilde E(t),\widetilde F(t))$ is a normal moving frame of
$\Lambda(t)$ if and only there exists a constant orthogonal matrix
$U$ of size $n\times n$ such that
\begin{equation}\label{3U}
\widetilde E(t)=E(t)U,\ \widetilde F(t)=F(t)U.
\end{equation}
\end{theor}
\begin{remark}
 \label{lcrem}
 Note that from \eqref{3structeq} it follows that if $\bigl(\widetilde E(t), \widetilde F(t)\bigr)$ is a Darboux moving frame such that
 $\widetilde E(t)$ is an orthonormal frame of $\Lambda(t)$  and ${\rm span}\,\{\widetilde F(t)\}=\Lambda^{\rm trans}(t)$. Then there exists a curve of antisymmetric matrices $B(t)$ such that
\begin{equation}
 \label{3structRiem1}
\left\{\begin{array}{l}
 \widetilde E'(t)=\widetilde E(t) B(t)+\widetilde F_a(t)\\
 \widetilde F'(t)=-\widetilde E(t)\widetilde{\mathcal R}(t)+\widetilde F(t) B(t),
\end{array}\right.
\end{equation}
where $\widetilde{\mathcal R}(t)$ is the matrix of the curvature map
$\mathfrak R(t)$ on $\Lambda(t)$ w.r.t. the basis $\widetilde E(t)$.
\end{remark}

As a matter of fact, normal moving frames define a principal
$O(n)$-bundle of symplectic frame in $G$ endowed with a canonical
connection. Also, relations \eqref{3U} imply that the following
$n$-dimensional subspaces
\begin{equation}\label{3V}
\Lambda^{\rm trans}(t)={\rm span}\{F(t)\}
\end{equation}
of $G$ does not depend on the choice of the normal moving frame. It
is called the \emph{canonical complement} of $\Lambda(t)$ in $G$.
Moreover, the subspaces $\Lambda(t)$ and $\Lambda^{\rm trans}(t)$
are endowed with the \emph {canonical Euclidean structure} such that
the tuple of vectors $E(t)$ and $F(t)$ constitute an orthonormal
frame w.r.t. to it, respectively.

Finally, the linear map from $\Lambda(t)$ to $\Lambda(t)$ with the
matrix $R(t)$ from (\ref{3structeq}) in the basis $\{E(t)\}$, is
independent of the choice of normal moving frames. It will be
denoted by $\mathfrak R(t)$ and it is called the \emph {curvature
map} of the curve $\Lambda(t)$.

\medskip
Now we apply the above results for curves in Lagrange Grassmannians
to sub-Riemannian structures. Since $\mathfrak J^{c}_\lambda(0)$ and
$\Pi^{c}_\lambda$ can be naturally identified, there is a canonical
splitting of $W^{c}_{\lambda}:$
  \begin{equation}\label{3splitting}
  W^{c}_{\lambda}=\Pi^{c}_\lambda\oplus\widetilde{\mathfrak J}^{c}(\lambda),
  \end{equation}
   where 
   $\widetilde {\mathfrak J}^{c}(\lambda)={\rm span}(F^{\lambda}(0))$ is the canonical complement.
In other words, $\widetilde{\mathfrak J}^{c}(\lambda)$ is actually a
(nonlinear) Ehresmann connection of $\Pi^{c}_\lambda$ in
$W^{c}_{\lambda}$. It also follows that the subspaces
$\Pi_\lambda^c$ and $\widetilde{\mathfrak J}^{c}(\lambda)$ are
equipped with a canonical Euclidean structure. Moreover, one can
define the curvature map of the dynamical Lagrangian distribution
$(\vec h,\Pi^c)$, i.e. $\mathfrak R^{c}_\lambda:
\Pi^{c}_\lambda\rightarrow \Pi^{c}_\lambda$ such that $\mathfrak
R^{c}_\lambda=\mathfrak R_\lambda(0)$ of the curvature maps of the
Jacobi curve $\mathfrak J^{c}_\lambda(\cdot)$ at $t=0$.
 This curvature maps are intrinsically related to the sub-Riemannian structure and will be called \emph{the reduced curvature map}
  of the sub-Riemannian structure.

\medskip
Let $\lambda\in T^*M$ and let $\lambda(t)=e^{t\vec h}\lambda$.
Assume that $(E^\lambda(t),F^\lambda(t))$ is a normal moving frame
of the Jacobi curve $\mathfrak J^{c}_\lambda(t)$ attached at point
$\lambda$. Let $\mathfrak E$ be the Euler field on $T^*M$, i.e. the
infinitesimal generator of the homotheties  on its fibers. Clearly
$T_\lambda(T^*M)=T_\lambda\mathcal H_{c}\oplus\mathbb R \mathfrak
E(\lambda)\oplus{\rm span}\{\partial_{u_i}(\lambda), 1\leq i\leq
s\}$. The flow $e^{t\vec h}$ on $T^*M$ induces the push-forward maps
$e^{t\vec h}_*$ between the corresponding tangent spaces $T_\lambda
T^*M$ and $T_{\lambda(t)}T^*M$, which in turn induce naturally the
maps between the spaces $T_\lambda(T^*M)/{\rm span}\{\vec
h(\lambda), \vec u_i(\lambda), 1\leq i\leq s\}$ and
$T_{\lambda(t)}T^*M/{\rm span}\{\vec h(\lambda(t)), \vec
u_i(\lambda(t)), 1\leq i\leq s\}$. The map $\mathcal K^t$ between
$T_\lambda(T^*M)/{\rm span}\{\vec h(\lambda), \vec u_i(\lambda),
1\leq i\leq s\}$ and $T_{\lambda(t)}T^*M/{\rm span}\{\vec
h(\lambda(t)), \vec u_i(\lambda(t)), 1\leq i\leq s\}$, sending
$E^\lambda(0)$ to $e^{t\vec h}_*E^{\lambda}(t)$, $F^\lambda(0)$ to
$e^{t\vec h}_*F^{\lambda}(t)$, and the equivalence class of
$\mathfrak E(\lambda), \partial_{u_i}(\lambda) (1\leq i\leq s)$ to
the equivalence class of $\mathfrak E(e^{t\vec h}\lambda),
\partial_{u_i}(\lambda(t)) (1\leq i\leq s)$, is independent of the
choice of normal moving frames. The map $\mathcal K^t$ is called
\emph{the parallel transport} along the extremal $e^{t\vec
h}\lambda$ at time $t$.  For any $v\in T_\lambda(T^*M)/{\rm
span}\{\vec h(\lambda), \vec u_i(\lambda), 1\leq i\leq s\}$, its
image $v(t)=\mathcal K^t(v)$ is called \emph{ the parallel transport
of $v$ at time $t$}. Note that from the definition of the reduced
Jacobi curves and the construction  of normal moving frames it
follows that the restriction of the parallel transport $\mathcal
K_t$ to the vertical subspace $T_\lambda(T_{\pi(\lambda)}^*M)$ of
$T_\lambda(T^*M)$ can be considered as a map onto the vertical
subspace $T_{\lambda(t)}(T_{\pi(\lambda(t))}^*M)$ of
$T_{\lambda(t)}(T^*M)$. A vertical vector field $V$ is called
\emph{parallel} if $V(e^{t\vec h}\lambda)=\mathcal
K^t\bigl(V(\lambda)\bigr)$.

%

{\bf Example 6} (Riemannian geodesic flow) In this case,
$\mathcal{D}=TM, W=0$ and there is no symmetries at all ($s=0$).
 In \cite{agfeedback} the reduced curvature map was expressed by the Riemannian curvature tensor. 
If we adopt the notations from Example 3 and take the constants
$c_0=\frac{1}{2}, c_i=0 (1\leq i\leq s)$. Then
\begin{equation}\label{Rie}
\mathfrak R_{\lambda}^c(v)=R^\nabla(p^h, v^h)p^h,\quad\forall
\lambda=(q, p)\in \mathcal H_{c} ,q\in M, p\in T^*_qM, v \in
\Pi^c_\lambda.
\end{equation}
Given a vector $X\in T_qM$ denote by $\nabla_X$ its lift to the
Levi-Civita connection, considered as an Ehresmann connection on
$T^*M$. Then by constructions the Hamiltonian vector field $\vec  h$
is horizontal and satisfies $\vec  h=\nabla_{p^h}$. Take any $v, w
\in \Pi_\lambda^c$ and let $V$ be a vertical vector field such that
$V(\lambda)=v$. From \eqref{Rie} , structure equation
\eqref{3structeq}, and the fact that the Levi-Civita connection (as
an Eheresmann connection on $T^*M$) is a Lagrangian distribution
(c.f. \cite{agfeedback}) it follows that the Riemannian curvature
tensor satisfies the following identity:
\begin{equation}
\label{Riemnabla} \langle R^\nabla(p^h, v^h)p^h,
w^h\rangle=-\sigma\left([\nabla_{p^h},\nabla_{V^h}](\lambda),\nabla
_{w^h}\right).\quad\square
\end{equation}
%

\subsection{Proof of Theorem \ref{3main1}}
We first express the canonical complement in terms of the
Levi-Civita connection of the Riemannian metric and the tensor
$J_i^{c}$ and then we can give the proof of Theorem \ref{3main1}
using some calclus formulae which is developped in \cite{cijacobi}.
\subsubsection{The canonical complement $\widetilde
{\mathfrak J}^{c}(\lambda)$} The restriction of the parallel
transport $\mathcal K^t$ to $\Pi^{c}_\lambda$ is characterized by
the following two properties:
\begin{enumerate}
\item $\mathcal K^t$ is an orthogonal transformation of spaces $\Pi_{\lambda}^{c}$ and $\Pi_{e^{t\vec h}\lambda}^{c}$;
\item The space ${\rm span} \{\frac{d}{dt}\bigl((e^{-t\vec  h})_*(\mathcal K^t v)\bigr)|_{_{t=0}}: v\in \Pi_{\lambda}^{c}\}$ is isotropic.
\end{enumerate}
Then $\widetilde{\mathfrak J}^{c}(\lambda)= {\rm span}
\{\frac{d}{dt}\bigl((e^{-t\vec  h})_*(\mathcal K^t
v)\bigr)|_{_{t=0}}: v\in \Pi_{\lambda}^{c}\}$.

To express $\widetilde {\mathfrak J}^{c}(\lambda)$ in terms of the
Riemannian manifold and the magnetic field, we show the
decomposition of the symplectic form $\sigma$ (the standard
symplectic form on $T^*\widetilde M$) and the Hamiltonian field
$\vec h$. One can see that the diffeomorphism $I^{c}$, defined as
before, are not in general symplectic. Indeed, each level set
inherits a symplectic structure depending on the choice of the level
 $\{c_i:1\leq i\leq s\}$.

By the construction of the map $I^{c}$, for any vector field $X$ on
$T^*\widetilde M$, we can assign the vector field $\underline{X}$ on
$T^*M$ s.t. $PR^{c}_*\underline X=((I^{c})^{-1})_*X$ and
$\pi_*\underline X\in\mathcal{D}$. In the following, denote by
$\nabla_{p^h}$ the lift of $p^h$ to $T^*\widetilde M$ with respect
to the Levi-Civita connection and denote
$\Omega^c=\sum_{i=1}^sc_id\omega_i$ and
$\bar\sigma=(I^{c}\circ\rm{PR}^{c})^*\sigma$. We will denote by
$\tilde\sigma$ the standard symplectic form on $T^*\widetilde M$.
The proof of the following lemma is complete similar to that of
Lemma 3.1-3.3 in \cite{cijacobi} and thus is omitted.

\begin{lemma}\label{3symplectic}
The following decomposition formulae hold.
\begin{enumerate}
\item
 On the level set
$\{u_i=c_i, i=1,...,s\},\quad \sigma=\bar\sigma-(\pi\circ{\rm
 pr})^*(\Omega^c)$;
\item For any  vectors $X, V \in T_\lambda T^*M$ with
$\pi_*V =0$ we have
$\sigma(X,v)=g(\pi_*X, V^h);$
\item $\vec  h(p,q)=\underline{\nabla_{p^h}}-(J^cp^h)^v+\overrightarrow W$.
\end{enumerate}
\end{lemma}

Comparing with the sub-Riemannian geodesic problem, we will develop
some additional calculus formulae about the potential $W$.
\begin{lemma}\label{CalW}
Let $V_1,V_2$ be the vector fields on $T^*M$ with
$\pi_*V_1=\pi_*V_2=0$. Then
\begin{enumerate}
\item $\overrightarrow W=-(\nabla W)^v;$
\item $\bar\sigma([\overrightarrow
W,\underline{\nabla_{V_1^h}}],\underline{\nabla_{V_2^h}}]=-{\rm
Hess}\ W(V_1^h,V_2^h);$
\item $\overrightarrow W\left(g(V_1^h,V_2^h)\right)-g\left(([\overrightarrow W,V_1])^h,V_2^h\right)-g\left(V_1^h,([\overrightarrow
W,V_2])^h\right)=0.$
\end{enumerate}
\end{lemma}
\begin{proof}  (1) From the
second item of the last lemma it follows that $$
\sigma(\overrightarrow W,V_1^h)=V_1^h(W)=dW(V_1^h)=g(\nabla
W,V_1^h)=-\sigma((\nabla W)^v, V_1^h).$$ Taking into account that
$\pi_*\overrightarrow W=0$, we get $$ \overrightarrow W=-(\nabla
W)^v,\quad{\rm span}\{\partial_{u_i},i=1,...s\}.$$

On the other hand, it follows from \eqref{simW} and item (1) of the
present lemma that
$$\sigma(\vec u_i,\overrightarrow W)=-\vec u_i(W)=-X_i(W)=0,\ i=1,...,s.$$
Thus, we get the required identity $\overrightarrow W=-(\nabla
W)^v$.

 (2) Both sides of the required identity are linear w.r.t.
$V_1,V_2$, respectively, thus it is sufficient to prove it for the
case that $V_1^h,V_2^h$ are both vector fields on $T^*\widetilde M$.
But for this case the required identity is a direct consequence of
the definition of the Hessian.

(3) Left-hand side is linear  w.r.t. $V_1,V_2$, respectively, thus
it is sufficient to prove it for the case that $V_1^h,V_2^h$ are
both vector fields on $T^*\widetilde M$. In this case, the vector
fields $(V_1^h)^v,(V_2^h)^v$, together with $(\nabla W)^v$ are all
constants on the fibers of $T^*M$ and then the required identity
become trivial.
\end{proof}
Given any $X\in \Pi_{\lambda}^{c}$ denote by $\widetilde
\nabla_{X^h}$ the lift of $X$ to $\widetilde {\mathfrak
J}^{c}(\lambda)$: the unique vector $\widetilde\nabla_{X^h}\in
\widetilde {\mathfrak J}^{c}(\lambda)$ such that $({\rm pr}\circ
\pi)_*\widetilde\nabla_{X^h}=X^h$. Then there exist the unique $B\in
{\rm End}(\Pi^{c}_{\lambda})$ and $\tilde A\in (\Pi_\lambda^c)^*$
such that
\begin{equation}
\label{decompc} \widetilde\nabla_{v^h}=\underline{\nabla_{v^h}}+Bv
,\quad \forall v\in \Pi_{\lambda}^{c},
\end{equation}
where $\nabla$ stands for the lifts to the Levi-Civita connection on
$T^*\widetilde M$, as before.
\begin{lemma}\label{anti}
The linear operator $B$ is antisymmetric w.r.t. the canonical
Euclidean structure in $\Pi_{\lambda}^{c}$.
\end{lemma}

\begin{proof}
Fix a point $\bar\lambda\in T^*M$ and consider a small neighborhood
$U$ of $\bar\lambda$. Let $\mathcal E=\{\mathcal E^i\}_{i=1}^{m-1}$
be a frame of $\Pi_{\lambda}^{c}$ (i.e. $\Pi_{\lambda}^{c}={\rm
span}\, \mathcal E(\lambda)$) for any $\lambda\in U$ such that the
following four conditions hold
\begin{enumerate}
\item $\mathcal E$ is orthogonal w.r.t. the canonical Euclidean structure on $\Pi_{\lambda}^{c}$;
\item Each vector field $\mathcal E^i$ is parallel w.r.t the canonical parallel transport $\mathcal K_t$, i.e. $\mathcal E^i(e^t\vec  h\lambda)=\mathcal K^t\mathcal E^i(\lambda)$ for any $\lambda$ and $t$ such that $\lambda, e^{t\vec  h}\lambda\in U$;
\item The vector fields $(J^cp^h)^v$ and $\mathcal E^i$ commute on $U\cap T_{\pi(\bar\lambda)}^*M$;
\item The vector fields $\vec  u_i,\ \forall 1\leq i\leq s$ and $\mathcal E^i$ commute on $U\cap T_{\pi(\bar\lambda)}^*M$.
\end{enumerate}
Note that the frame $\mathcal E$ with properties above exists,
because the Hamiltonian vector field $\vec  h$ is transversal to the
fibers of $T^*M$ and it commutes with $\vec u_i,\ \forall 1\leq
i\leq s$.

From the property (2) of the parallel transport $\mathcal K^t$ in
this subsection it follows that
\begin{equation}\label{3trans2}
\widetilde\nabla_{(\mathcal E^i)^h}=-{\rm ad}\vec  h\, \mathcal E^i
\end{equation}

Using the above defined identification $I^{c}:N^c\to T^*\widetilde
M$, one can look on the restriction of the tuple of vector fields
$\mathcal E$ to the submanifold $\{u_i=c_i, i=1,...,s\}$ as on the
tuple of the vertical vector fields of $T^*\widetilde M$ (which
actually span the tangent to the intersection of the fiber of
$T^*\widetilde M$ with the level to the corresponding Riemannian
Hamiltonian). Then first the tuple $\mathcal E$ is the tuple of
orthonormal vector fields (w.r.t. the canonical Euclidean structure
on the fibers of $T^*\widetilde M$, induced by the Riemannian metric
$g$). Further, by the equations \eqref{3structRiem1} the Levi-Civita
connection of $g$ is characterized by the fact that there exists a
field of antisymmetric operators $\widetilde B\in {\rm
End}(\Pi^{c}_\lambda)$ such that
\begin{equation}
\label{3nablA} [\nabla_{p^h},\widetilde{\mathcal
E}^i(\lambda)]=-\nabla_{\bigl(\widetilde{\mathcal
E}^i(\lambda)\bigr)^h} -\widetilde B \widetilde{\mathcal
E}^i(\lambda)
\end{equation}
On the other hand, from \eqref{3trans2},\eqref{3nablA}, using the
second item of Lemma \ref{3symplectic} and the property (3) of
$\mathcal E^i$, one has
\begin{equation}
\label{3seria1}
\begin{split}
\widetilde\nabla_{(\mathcal E^i)^h}=-{\rm ad}\vec  h\, \mathcal E^i=
-\bigl[\underline{\nabla_{p^h}}-\sum_{i=1}^sc_i(J_ip^h)^v+\overrightarrow
W,\mathcal E^i]=\underline{\nabla_{\bigl(\mathcal
E^i(\lambda)\bigr)^h}} +\widetilde B\,\mathcal
E^i(\lambda)+[\overrightarrow W,\mathcal E^i(\lambda)].
\end{split}
\end{equation}
From item (3) of Lemma \ref{CalW} and property (1) of $\mathcal E^i$
it follows
\begin{equation}\label{OW}
g(([\overrightarrow W,\mathcal E^i])^h,(\mathcal E^j)^h)+
g((\mathcal E^i)^h,([\overrightarrow W,\mathcal
E^j])^h)=\overrightarrow W\left(g((\mathcal E^i)^h,(\mathcal
E^j)^h)\right)=0.
\end{equation}
Therefore, from \eqref{3seria1} and \eqref{OW} we conclude that $B$
is antisymmetric.
\end{proof}

\begin{lemma}\label{3liftB}
The operator $B$ satisfies
\begin{equation}
\label{3B1} (Bv)^h=-\frac{1}{2}J^cv^h\quad{\rm mod}\ p^h, \quad
\forall v\in \Pi^{c}_\lambda.
\end{equation}
\end{lemma}
\begin{proof}
Since $\widetilde {\mathfrak J}^{c}(\lambda)$ is an isotropic
subspace, we have
\begin{equation*}
\sigma(\widetilde\nabla_{v_1^h}, \widetilde\nabla_{v_2^h})=0,\quad
\forall\, v_1,v_2\in \Pi^{c}_\lambda.
\end{equation*}

On the other hand, using Proposition \ref{3symplectic}, the fact
that the Levi-Civita connection (as an Ehresmann connection) is a
Lagrangian distribution in $T^*\widetilde M$ and Lemma
\ref{3symplectic}, we get
\begin{eqnarray*}
0=\sigma(\widetilde\nabla_{v_1^h}, \widetilde\nabla^c_{v_2^h})&=&
\Bigl((I^{c}\circ\hbox{PR}^{c})^*\tilde\sigma-({\rm
pr}\circ\pi)^*\Omega^c\Bigr)\Bigl(\underline{\nabla_{v^h_1}}+Bv_1,
\underline{\nabla_{v^h_2}}+Bv_2\Bigr)\\
&=&-\Omega^c(v_1^h,v_2^h)-g\big((Bv_1)^h, v_2^h\big)+g\big((Bv_2)^h,
v_1^h)\\ &=&-g(J^cv_1^h, v_2^h)-g\big((Bv_1)^h,
v_2^h\big)+g\big((B^*v_1)^h, v_2^h).
\end{eqnarray*}
where $B^*$ is the dual of $B$ w.r.t. the Euclidean structure in
$\Pi_\lambda^c$. Taking into account that $B$ is antisymmetric, we
get
 \begin{equation}\label{3B2}
 (Bv)^h=-\frac{1}{2}J^cv^h\quad{\rm mod}\ p^h.
 \end{equation}
\end{proof}

\begin{cor}\label{3F}
The canonical complement $\widetilde {\mathfrak J}^{c}(\lambda)$ can
be expressed as follows:
$$\widetilde {\mathfrak J}^{c}(\lambda)=\{\underline{\nabla_{v^h}}-\frac{1}{2}(J^cv^h)^v-\frac{1}{2\|p^h\|^2}\cdot g\left(v^h,J^cp^h+2\nabla W\right)(p^h)^v,\ v\in
\Pi^{c}_\lambda\}.$$
\end{cor}
\begin{proof}
It follows from \eqref{decompc} and \eqref{3B1} that there exist
$A\in (\Pi^c_\lambda)^*$ such that
$$\widetilde\nabla_{v^h}=\underline{\nabla_{v^h}}-\frac{1}{2}(J^cv^h)^v+A(v)(p^h)^v.$$
Note that $\sigma(\vec h,(p^h)^v)=g(p^h,p^h)=\|p^h\|^2$. Hence, from
the fact that $\widetilde\nabla_{v^h}$ is tangent to the Hamiltonian
vector field $\vec h$, we get easily that
$$A(v)=-\frac{1}{2\|p^h\|^2}\cdot g\left(v^h,J^cp^h+2\nabla W\right),$$which completes the proof of the corollary.
\end{proof}

\subsubsection{The reduced curvature map}

As a direct consequence of structure equation \eqref{3structeq}, we get the following preliminary descriptions of the reduced curvature map:
\begin{prop}\label{preli}
Let $v\in \Pi^{c}_\lambda$. Let $V$ be a parallel vector field such
that $V(\lambda)=v$. Then the curvature maps satisfy the following
identities:
\begin{eqnarray}\label{3repreli}
g\big((\mathfrak R_{\lambda}^{c}v)^h,v^h\big)=-\sigma(\hbox{ad}\vec
h\ (\widetilde\nabla_{V^h}),\widetilde\nabla_{v^h}).
\end{eqnarray}
\end{prop}

It follows that in order to calculate the reduced curvature map it
is sufficient to know how to express the Lie bracket of vector
fields on the cotangent bundle $T^*M$ via the covariant derivatives
of Levi-Civita connection on $T^*\widetilde M$. 

\begin{prop}\label{3basic}
For any tensors $A, B$ of type $(1,1)$ on $\widetilde M$, the
following identity holds:
\begin{enumerate}
\item $[(Ap^h)^v, (Bp^h)^v]=(B(Ap^h))^v-(A(Bp^h))^v,$
\item $[\underline{\nabla_{p^h}}, (Ap^h)^v]=-\underline{\nabla_{Ap^h}}+((\nabla_{p^h} A)p^h)^v.$
\end{enumerate}
\end{prop}
For simplicity, denote $\bar\sigma=(I^{c}\circ\rm{PR}^{c})^*\sigma$
and $\Omega^c=\sum_{i=1}^sc_id\omega_i,$ as before.
As in the proof of Lemma \ref{anti}, we can take a parallel vector
field $V$ such that $V(\lambda)=v$ and
\begin{equation}\label{3commute}
[(J^{c}p^h)^v, V](\bar\lambda)=0,\quad\bar\lambda\in U\cap T_q^*M,
\end{equation}
where $U$ is a neighborhood of $\lambda.$ Similar to Proposition 4
of \cite{cijacobi}, we have

\begin{lemma}
\label{3Riesimp} Let $V,V_1,V_2$ be vector fields on $T^*M$ with
$\pi_*V=\pi_*V_1=\pi_*V_2=0$. Then
\begin{enumerate}
\item $([(J^{c}p^h)^v, (J^{c}V^h)^v])^h=J^c([(J^{c}p^h)^v, (V^h)^v])^h$,
\item $\bar\sigma([(J^{c}p^h)^v, \underline{\nabla_{V_1^h}}],\underline{\nabla_{V_2^h}})=g(\nabla J^{c}(p^h,V_1^h),V_2^h),$
\item $({\rm pr}\circ\pi)_*([(J^{c}p^h)^v, \underline{\nabla_{V^h}}])=J^cV^h,$
\item $({\rm pr}\circ\pi)_*([\underline{\nabla_{p^h}}, \underline{\nabla_{V^h}}])=\frac{1}{2}J^{c}V^h-\frac{1}{2\|p^h\|^2}g(J^{c}V^h,p^h)p^h$.
\end{enumerate}
\end{lemma}

Let us simplify the right-hand side of the identity
\eqref{3repreli}. First, from the last line of the structural
equations \eqref{3structeq} it follows that
\begin{equation}\label{3simp1}
({\rm pr}\circ\pi)_*(\hbox{ad}\vec  h(\widetilde\nabla_{V^h}))\in
\mathbb R p^h.
\end{equation}
Besides,
\begin{equation*}
\sigma\left(\vec
h,\frac{1}{2}(J^cv^h)^v+\frac{1}{2\|p^h\|^2}g\left(v^h,J^cp^h\right)(p^h)^v\right)
=\frac{1}{2}g(p^h,J^cv^h)+\frac{1}{2\|p^h\|^2}g(v^h,J^cp^h)g(p^h,p^h)=0.
\end{equation*}
 Hence from Lemma \ref{3symplectic} and Corollary \ref{3F} it
follows that
\begin{eqnarray*}
 &&\sigma(\hbox{ad}\vec  h(\nabla_{V^h}),
\nabla_{v^h})= \sigma\left(\hbox{ad}\vec  h(\nabla_{V^h}),
\underline{\nabla_{v^h}}-\frac{1}{\|p^h\|^2}g(v^h,\nabla
W)(p^h)^v\right)\\
&=&\sigma\left([\underline{\nabla_{p^h}}-(J^cp^h)^v,\underline{\nabla_{V^h}}-\frac{1}{2}(J^cV^h)^v-\frac{1}{2\|p^h\|^2}\cdot
g(V^h,J^cp^h)(p^h)^v],\underline{\nabla_{v^h}}\right)\\
&+&\sigma\left([\underline{\nabla_{p^h}}-(J^cp^h)^v,\underline{\nabla_{V^h}}-\frac{1}{2}(J^cV^h)^v-\frac{1}{2\|p^h\|^2}\cdot
g(V^h,J^cp^h)(p^h)^v],-\frac{1}{\|p^h\|^2}\cdot g(v^h,\nabla
W)(p^h)^v\right)\\
&+&\sigma\left([\underline{\nabla_{p^h}}-(J^cp^h)^v,-\frac{1}{\|p^h\|^2}\cdot
g(V^h,\nabla
W)(p^h)^v],\underline{\nabla_{v^h}}-\frac{1}{\|p^h\|^2}g(v^h,\nabla
W)(p^h)^v\right)\\
 &+&\sigma\left([\overrightarrow
W,\underline{\nabla_{V^h}}],\underline{\nabla_{v^h}}-\frac{1}{\|p^h\|^2}g(v^h,\nabla
W)(p^h)^v\right)\\
&+&\sigma\left([\overrightarrow
W,-\frac{1}{2}(J^cV^h)^v-\frac{1}{2\|p^h\|^2}\cdot
g(V^h,J^cp^h+2\nabla
W)(p^h)^v],\underline{\nabla_{v^h}}\right)\\
&=:&T_1+T_2+T_3+T_4+T_5
\end{eqnarray*}

\medskip
We will deal with the terms $T_i (1\leq i\leq 5)$ in steps.

\noindent\textit{Step 1} It follows identity \eqref{Riemnabla} that
\begin{equation}
\label{3simp4} \bar\sigma([\underline{\nabla_{p^h}},
\underline{\nabla_{V^h}}],\underline{\nabla_{v^h}})=-g(R^\nabla(p^h,v^h)p^h,v^h).
\end{equation}
Also it follows from item (3) of Lemma \ref{3Riesimp} and item (4)
of Lemma \ref{3Riesimp} that
\begin{eqnarray*}
&&\Omega^{c}(({\rm pr}\circ\pi)_*([\underline{\nabla_{p^h}},
\underline{\nabla_{V^h}}]),
v^h)+\frac{1}{2}\bar\sigma([\underline{\nabla_{p^h}}, (J^cV^h)^v],
\underline{\nabla_{v^h}})\\
&=&-g\bigl(({\rm pr}\circ\pi)_*([\underline{\nabla_{p^h}},
\underline{\nabla_{V^h}}]), J^{c}v^h\bigr)+\frac{1}{2}\bar\sigma([\underline{\nabla_{p^h}},\underline{\nabla_{V^h}}],(J^{c}v^h)^v)\\
&=&-\frac{1}{2}g\bigl(({\rm
pr}\circ\pi)_*([\underline{\nabla_{p^h}},
\underline{\nabla_{V^h}}]), J^{c}v^h\bigr)\\
&=&-\frac{1}{4}\|J^cv^h\|^2+\frac{1}{4\|p^h\|^2}\left(g(J^cv^h,p^h)\right)^2
\end{eqnarray*}
Also it follows from straightforward computations that
\begin{equation}
\label{3simp11} \Omega^{c}\bigl(({\rm
 pr}\circ\pi)_*([\underline{\nabla_{p^h}},
(J^{c}v^h)^v]),v^h\bigr)=-\Omega^{c}\bigl(J^{c}v^h,v^h\bigr)=\|J^{c}v^h\|^2
\end{equation}
Also it follows from item (2) of Proposition \ref{3basic} that
\begin{equation}
\label{3simp6} \bar\sigma([\underline{\nabla_{p^h}},
(p^h)^v],\underline{\nabla_{v^h}})=\bar\sigma(-\underline{\nabla_{p^h}},\underline{\nabla_{v^h}})=0.
\end{equation}
and
\begin{equation}
\label{3simp12} \Omega^{c}\bigl(({\rm
pr}\circ\pi)_*([\underline{\nabla_{p^h}},
(p^h)^v]),v^h\bigr)=-\Omega^{c}(p^h,v^h\bigr) =g(p^h, J^cv^h).
\end{equation}
It follows from item (2) of Lemma \ref{3Riesimp} that
 \begin{equation}
\label{3simp5} \bar\sigma([(J^{c}p^h)^v,
\underline{\nabla_{V^h}}],\underline{\nabla_{v^h}})=g(\nabla
J^{c}(p^h,v^h),v^h).
\end{equation}
Applying item (3) of Lemma \ref{3Riesimp}, we get
\begin{equation}
\label{3simp10} \Omega^{c}(({\rm pr}\circ\pi)_*([(J^{c}p^h)^v,
\underline{\nabla_{V^h}}]),v^h)=\Omega^{c}(J^{c}v^h,v^h)=-\|J^{c}v^h\|^2.
\end{equation}
And it follows from \eqref{3commute} and item (1) of Lemma
\ref{3Riesimp} that
\begin{equation}\label{simW}
[(J^cp^h)^v,(J^cV^h)^v]=0.
\end{equation}
And it follows from item (1) of Proposition \ref{3basic} that
\begin{equation}\label{Jcp}
[(J^cp^h)^v,(p^h)^v]=0.
\end{equation}
Then it follows that
\begin{equation*}
\sigma\left([-(J^cp^h)^v,-\frac{1}{2\|p^h\|^2}\cdot
g(V^h,J^cp^h+2\nabla W)(p^h)^v],\underline{\nabla_{v^h}}\right)=0
\end{equation*}
Summarizing all the calculations above, we have
\begin{equation}\label{T1}
T_1=-g(R^\nabla(p^h,v^h)p^h,v^h)-g(\nabla
J^{c}(p^h,v^h),v^h)-\frac{1}{4}\|J^cv^h\|^2-\frac{3}{4\|p^h\|^2}\left(g(p^h,J^cv^h)\right)^2
\end{equation}
\medskip
\textit{Step 2} Again, it follows from item (4) of  Lemma
\ref{3Riesimp} that
\begin{equation}
\sigma([\underline{\nabla_{p^h}},\underline{\nabla_{V^h}}],(p^h)^v)=-\frac{1}{2}g(J^{c}v^h,p^h).
\end{equation}
And it follows from straightforward computations that
\begin{equation}
\sigma\left([\underline{\nabla_{p^h}},\frac{1}{2}(J^cV^h)^v+\frac{1}{2\|p^h\|^2}\cdot
g(V^h,J^cp^h)(p^h)^v],(p^h)^v\right)=0.
\end{equation}
And it follows from item (3) of Lemma \ref{3Riesimp} that
\begin{equation}
\sigma\left([(J^cV^h)^v,\underline{\nabla_{V^h}}],(p^h)^v\right)=g(J^cv^h,p^h)).
\end{equation}
 Hence, we have
\begin{equation}\label{T2}
T_2=-\frac{1}{\|p^h\|^2}g(J^cp^h,v^h)g(\nabla W,v^h).
\end{equation}
\medskip
\textit{Step 3} It follows from item (2) of Proposition \ref{3basic}
that $[\underline{\nabla_{p^h}},(p^h)^v]=-\underline{\nabla_{p^h}}$,
hence
\begin{eqnarray*}
&&\sigma\left([\underline{\nabla_{p^h}},-\frac{1}{\|p^h\|^2}\cdot
g(V^h,\nabla
W)(p^h)^v],\underline{\nabla_{v^h}}-\frac{1}{\|p^h\|^2}g(v^h,\nabla
W)(p^h)^v\right)\\
&=&-\frac{1}{\|p^h\|^2}\cdot g(v^h,\nabla
W)\sigma\left([\underline{\nabla_{p^h}},(p^h)^v],\underline{\nabla_{v^h}}-\frac{1}{\|p^h\|^2}g(v^h,\nabla
W)(p^h)^v\right)\\
&=&-\frac{1}{\|p^h\|^2}\cdot g(v^h,\nabla
W)\sigma\left(-\underline{\nabla_{p^h}},\underline{\nabla_{v^h}}-\frac{1}{\|p^h\|^2}g(v^h,\nabla
W)(p^h)^v\right)\\
&=&-\frac{1}{\|p^h\|^2}\cdot g(v^h,\nabla
W)g(J^cp^h,v^h)-\frac{1}{\|p^h\|^2}\cdot \left(g(V^h,\nabla
W)\right)^2.
\end{eqnarray*}
And it follows from \eqref{Jcp} that
\begin{equation}
\sigma\left([(J^cp^h)^v,-\frac{1}{\|p^h\|^2}\cdot g(V^h,\nabla
W)(p^h)^v],\underline{\nabla_{v^h}}\right)=0.
\end{equation}
Hence,
\begin{equation}\label{T3}
T_3=-\frac{1}{\|p^h\|^2}\cdot g(v^h,\nabla
W)g(J^cp^h,v^h)-\frac{1}{\|p^h\|^2}\cdot \left(g(V^h,\nabla
W)\right)^2
\end{equation}
\medskip
\textit{Step 4} From item (2) of Lemma \ref{CalW}, we have
\begin{equation}
\bar\sigma\left([\overrightarrow
W,\underline{\nabla_{V^h}}],\underline{\nabla_{v^h}}\right)=-{\rm
Hess}\ W(v^h,v^h).
\end{equation}

\begin{lemma}\label{Pw}
The following identity holds.
$$({\rm pr}\circ\pi)_*([\overrightarrow W,
\underline{\nabla_{V^h}}])=\frac{1}{\|p^h\|^2}g(v^h,\nabla W)p^h.$$
\end{lemma}
\begin{proof}
First of all, from \eqref{3simp1},item (3) of Lemma
\ref{3symplectic} and Lemma \ref{3F} it follows
$$({\rm pr}\circ\pi)_*([\underline{\nabla_{p^h}}-(J^cp^h)^v+\overrightarrow W,
\underline{\nabla_{V^h}}-\frac{1}{2}(J^cV^h)^v])=0,\quad{\rm mod}\
p^h.$$ Then together with item (3)-(4) of Lemma \ref{3Riesimp} and
item (2) of Proposition \ref{3basic} we have
$$({\rm pr}\circ\pi)_*([\overrightarrow W,
\underline{\nabla_{V^h}}])=0,\quad{\rm mod}\ p^h.$$

Furthermore, from the classical Cartan's formula we have
\begin{eqnarray}\label{Cartan}
0&=&d\sigma(\overrightarrow
W,\underline{\nabla_{V^h}},(p^h)^v)=\overrightarrow
W(\sigma(\underline{\nabla_{V^h}},(p^h)^v))-\underline{\nabla_{V^h}}(\sigma(\overrightarrow
W,(p^h)^v))+(p^h)^v(\sigma(\overrightarrow
W,\underline{\nabla_{V^h}}))
\\
\notag &-&\sigma([\overrightarrow
W,\underline{\nabla_{V^h}}],(p^h)^v)+\sigma([\overrightarrow
W,(p^h)^v],\underline{\nabla_{V^h}})-\sigma([\underline{\nabla_{V^h}},(p^h)^v],\overrightarrow
W).
\end{eqnarray}
Since $\overrightarrow W$ is constant on the fiber of $T^*M$, one
can easily show $$(p^h)^v(\sigma(\overrightarrow
W,\underline{\nabla_{V^h}}))-\sigma([\underline{\nabla_{V^h}},(p^h)^v],\overrightarrow
W)=0.$$ Then, $$\sigma([\overrightarrow
W,\underline{\nabla_{V^h}}],(p^h)^v)=\sigma([\overrightarrow
W,(p^h)^v],\underline{\nabla_{V^h}})=g(\nabla W,v^h).$$ Hence, the
required identity follows and the lemma is proved.
\end{proof}
As a direct consequence of the last lemma, we have
\begin{eqnarray*}
\Omega^{c}\bigl(({\rm
 pr}\circ\pi)_*([\overrightarrow W,
\underline{\nabla_{V^h}}]),v^h\bigr)&=&\frac{1}{\|p^h\|^2}g(J^cp^h,v^h)g(v^h,\nabla
W),\\
\sigma\left([\overrightarrow
W,\underline{\nabla_{V^h}}],(p^h)^v\right)&=&g(v^h,\nabla W).
\end{eqnarray*}
 As a result of above calculations, we get $$T_4=-{\rm Hess}\
W(v^h,v^h)-\frac{1}{\|p^h\|^2}g(J^cp^h,v^h)g(v^h,\nabla
W)-\frac{1}{\|p^h\|^2}(g(v^h,\nabla W))^2.$$
\medskip
\textit{Step 5} First of all, we show the following
\begin{lemma}\label{JPw}
The following identity holds.
$$[\overrightarrow W,(J^cV^h)^v]=\frac{1}{\|p^h\|^2}g(v^h,\nabla W)(Jp^h)^v.$$
\end{lemma}
\begin{proof}
As $\overrightarrow W=-(\nabla W)^v$ is constant on the fiber of
$T^*M$, we can proceed with the following calculations
\begin{equation}
([\overrightarrow W,(J^cV^h)^v])^h=({\rm
pr}\circ\pi)_*([\overrightarrow
W,\underline{\nabla_{J^cV^h}}])=J^c\left(({\rm
pr}\circ\pi)_*([\overrightarrow W,\underline{\nabla_{V^h}}])\right).
\end{equation}
Substituting the identity of Lemma \ref{Pw} into the last identity,
we get $$([\overrightarrow
W,(J^cV^h)^v])^h=\frac{1}{\|p^h\|^2}g(v^h,\nabla W)J^cp^h,$$ and
then the required identity follows.
\end{proof}
As a direct consequence, we have
$$\sigma\left([\overrightarrow
W,(J^cV^h)^v],\underline{\nabla_{v^h}}\right)=-\frac{1}{\|p^h\|^2}g(v^h,\nabla
W)g(J^cp^h,v^h).$$ Furthermore, since $[\overrightarrow
W,(p^h)^v]=\overrightarrow W=-(\nabla W)^v$,then
$$
\sigma\left([\overrightarrow
W,(p^h)^v],\underline{\nabla_{v^h}}\right)=g(\nabla W,v^h).
$$
Therefore,
$$T_5=-\frac{1}{\|p^h\|^2}(g(v^h,\nabla
W))^2.$$

\medskip
Combining the results of Step 1-5 and using the fact
$\|p^h\|^2=2(c_0+W)$, we get the expression of the reduced curvature
maps, as shown in Theorem \ref{3main1}.

\bibliography{myrefe}
\end{document}